\documentclass[12pt]{article}
= 0.0in \headheight = 0.0in \topmargin = 0.3in
\def\Dj{\hbox{D\kern-.73em\raise.30ex\hbox{-}
\raise-.30ex\hbox{}}}
\def\dj{\hbox{d\kern-.33em\raise.80ex\hbox{-}
\raise-.80ex\hbox{\kern-.40em}}}
\usepackage{epsfig}
\usepackage{amsmath,amsthm,amsfonts,amssymb,amscd,cite}
\allowdisplaybreaks
\usepackage{mathrsfs}
\def \la{\lambda}

\newtheorem{theorem}{Theorem}

\newtheorem{lemma}[theorem]{Lemma}

\theoremstyle{remark}

\begin{document}
\title{On   graphs having minimal fourth adjacency coefficient }
\author{Shi-Cai Gong\thanks{Corresponding author. E-mail addresses:
scgong@zafu.edu.cn (S. Gong), sunshaowei2009 @126.com(S. Sun).} \thanks{ Supported by Zhejiang Provincial Natural Science Foundation of
China (No. LY20A010005), and National Natural Science Foundation of
China (No. 11571315,11601006).}, Peng Zou, Li-Ping Zhang~~and Shao-Wei Sun\thanks{ Supported by National Natural Science Foundation of
China (No. 11901525).}
\\{\small \it  School of Science, Zhejiang University of Science and Technology, }\\{\small \it
Hangzhou, 310023, P. R. China}
   }
\date{}
\maketitle

\baselineskip=0.20in

\noindent {\bf Abstract. } Let $G$ be a graph with order $n$ and adjacency matrix $\mathbf{A}(G)$. The adjacency polynomial of $G$ is defined as
$\phi(G;\la) =det(\la\mathbf{I}-\mathbf{A}(G))=\sum_{i=0}^n\mathbf{a_i}(G)\la^{n-i}$. Hereafter,
 $\mathbf{a}_i(G)$ is called the $i$-th adjacency coefficient of  $G$. Denote by $\mathfrak{G}_{n,m}$ the set of all connected graphs having $n$ vertices and $m$ edges.  A   graph  $G$ is said $4$-Sachs minimal if $$\mathbf{a}_4(G)=min\{\mathbf{a}_4(H)|H\in \mathfrak{G}_{n,m}\}.$$ The  value  $min\{\mathbf{a}_4(H)|H\in \mathfrak{G}_{n,m}\}$ is called the minimal $4$-Sachs number in $\mathfrak{G}_{n,m}$, denoted by $\bar{\mathbf{a}}_4(\mathfrak{G}_{n,m})$.
 \vspace{2mm}

In this paper, we study the relationship between the value $\mathbf{a}_4(G)$ and its structural properties. Especially, we give a structural characterization on  $4$-Sachs minimal graphs, showing that each $4$-Sachs minimal graph  contains a difference graph as its spanning subgraph (see Theorem \ref{12}). Then, for $n\ge 4$ and $n-1\le m\le 2n-4$,
  we  determine all $4$-Sachs minimal graphs together with the corresponding minimal $4$-Sachs number $\bar{\mathbf{a}}_4(\mathfrak{G}_{n,m})$.  \vspace{2mm}

\vspace{3mm}



\baselineskip=0.30in

\section{Introduction}
Let $G = (V,E)$ be a simple undirected graph with $\nu(G):=|V|=n$ vertices and $\epsilon(G):=| E|=m$ edges. Then  $G$  is called an \emph{$(n,m)$-graph}. Denote by $\mathfrak{G}_{n,m}$ and $\mathfrak{B}_{n,m}$ the set of all connected $(n,m)$-graphs and all connected bipartite $(n,m)$-graphs, respectively. The adjacency polynomial $\phi(G;\la)$ of $G$ is defined as
$$\phi(G;\la)=det(\la \mathbf{I}-\mathbf A(G))=\sum_{i=0}^n\mathbf{a_i}(G)\la^{n-i}.$$
Henceforth, we refer 
$\mathbf{a_i}(G)$, short for $\mathbf{a_i}$, as the $i$-th adjacency coefficient of $G$.\vspace{2mm}

 Denote by $o(G)$ and $c(G)$ the number of  \emph{components} and \emph{cycles} contained in the graph $G$, respectively. The subgraph $H$   of $G$  is called an \emph{$i$-Sachs subgraph} (of $G$) if  the order of $H$ is $i$ and each component of $H$ is either a single edge or a cycle. For each $i$, due to Sachs \cite{cds}, the adjacency coefficients $\mathbf{a_i}(G)$ of a graph $G$ can be expressed in terms of all its $i$-Sachs subgraphs of $G$ by the following result.\vspace{2mm}

\begin{theorem} \label{1} \em{\cite[Theorem 1.3]{cds}} Let $G$ be a graph with order $n$ and adjacency polynomial $\phi(G;\la)=\sum_{i=0}^n\mathbf{a_i}(G)\la^{n-i}$. Then
$$\mathbf{a_i}(G)=\sum_H(-1)^{o(H)}2^{c(H)},$$ where the
summation is over all $i$-Sachs subgraphs $H$ contained in $G$.
\end{theorem}

 Let $\la_1,\la_2,\ldots,\la_n$  be all eigenvalues of $ A(G)$.  From Viette¡¯s formulas, we have $$\mathbf{a}_i= \sigma_i(\la_1,\la_2,\ldots,\la_{n})=\sum_{I\subseteq\{1,2,\ldots,n\},|I|=i}\prod_{j\in I}\la_j.$$   In particular,  $\mathbf{a}_0 =1$, $\mathbf{a}_1 =0$, $\mathbf{a}_2 =-\epsilon(G)$, the opposite of the cardinality of edges contained in $G$, and $\mathbf{a}_3$ equals the number of triangles contained in $G$ multiplied by the constant $-2$.\vspace{2mm}

 An \emph{$r$-matching} in the graph $G$ is
a subset with $r$ edges such that every vertex of $V(G)$ is incident
with at most one edge in it.
The  $r$-\emph{matching number}, denoted by $\mathbf{m}_r(G)$, is defined as the cardinality of
$r$-matchings contained in $G$. Denote by $\mathbf{q}(G)$ the number of  quadrangles of $G$. Applying Theorem \ref{1}, we have
$$\mathbf{a}_4(G)=\mathbf{m}_2(G)-2\mathbf{q}(G).\eqno{(1.1)}$$
  From Eq.(1.1), for a given $(n,m)$-graph $G$, the adjacency coefficient $\mathbf{a}_4(G)$ is related to its structural properties,  not a fixed value.
Therefore, for a given graph $G$, it is interesting  to investigate the relationship between the value $\mathbf{a}_4(G)$ and its structural properties.\vspace{2mm}


A graph $G = (V,E)$ is said to be a \emph{threshold} graph if there exists a threshold $t$ and a function $w : V(G) \rightarrow R$ such that $uv\in E(G)$  if and only if $w(u) +w(v) \ge t$.
The graph $G $ is said to be a \emph{difference} graph if there exists a threshold $t$ and a function $w : V(G) \rightarrow R$ such that $|w(v)| < t$ for all $v \in V$ and distinct vertices $u$ and $v$ are adjacent if and only if $|w(u) - w(v)| \ge t$. \vspace{2mm}

Threshold graphs have a beautiful structure and possess many important
mathematical properties such as being the extreme cases of certain graph
properties.  For instance in the class $\mathfrak{G}_{n,m}$ threshold graphs maximize the number of independent sets and minimize the number of $k$-matchings; see \cite{cr,kr}. Threshold graphs also have applications in many areas such as computer
science and psychology. For more information on threshold graphs, one can see the book \cite{mp} and the references
therein. Difference graphs are called Threshold bipartite graphs in \cite{mp} and chain graphs in \cite{y}. A threshold graph can be obtained from a difference graph by adding
all possible edges in one of the partite sets (on either side). Therefore, threshold graphs and difference graphs are closely interconnected.\vspace*{2mm}

In \cite{gs}, Gong and Sun refer the fourth adjacency coefficient $\mathbf{a}_4(G)$ of the bipartite graph $G$  as the {\it $4$-Sachs number} of $G$. In \cite{gs}, Gong and Sun studied the structural properties of bipartite graphs having minimal  $4$-Sachs number among all bipartite graphs of $\mathfrak{B}_{n,m}$.
 Moreover, for $n\ge 6$ and $n-1\le m\le 2n-4$, the unique bipartite graph having minimal $4$-Sachs number in $\mathfrak{B}_{n,m}$ is determined in \cite{gs}. For consistency, in this paper the graph $G$ is called   a $4$-Sachs minimal $(n,m)$-graph if
$$\mathbf{a}_4(G)=\min\{\mathbf{a}_4(H)|H\in \mathfrak{G}_{n,m}\}.$$

In this paper,  we will investigate the relationship between the fourth adjacency coefficient  and the structural properties of   a given graph.
Especially, we give a structure characterization on  $4$-Sachs minimal $(n,m)$-graphs,  showing that each $4$-Sachs minimal $(n,m)$-graph  contains a  difference graph as its spanning subgraph.
 In addition, we   determine all  $4$-Sachs minimal $(n,m)$-graphs for $n\ge 6$ and $n-1\le m\le 2n-4$ together with the corresponding minimal $4$-Sachs number $\bar{\mathbf{a}}_4(\mathfrak{G}_{n,m})$.\vspace{2mm}

The rest of the paper is organized as follows: In section 2, we  give some notation and
essential preliminary results. Then  we show that each $4$-Sachs minimal $(n,m)$-graph  contains a  difference graph as its spanning subgraph in section $3$. In section $4$, we determine all $4$-Sachs minimal graphs  in $\mathfrak{G}_{n,m}$ together with the corresponding minimal $4$-Sachs number   for $n\ge 6$ and $n-1\le m\le 2n-4$.

\section{Preliminaries}
In this section, we  introduce some concepts, notations and preliminary results.
Let $G=(V,E)$ be a  graph and $v\in V(G)$. We use $N_G(v)$ to define the \emph{neighbor set} of $v$ in $G$, and let $d_G(v)=|N_G(v)|$ denote the \emph{degree} of $v$. If there is no confusion, we  simply $d_G(v)$ and $N_G(v)$ as $d(v)$ and $N(v)$, respectively. Denote by $V(e)$ the end-vertices of the edge $e$, {\it i.e.,} $V(e)=\{u,v\}$ if $e=uv$. As usually, the \emph{maximum} and \emph{minimum degree} of $G$ are written as $\Delta(G)$ and $\delta(G)$, respectively. Two vertices $u$ and $v$ of $G$ are called \emph{duplicate} if $N_G(u)=N_G(v)$. A vertex $v$ is   a \emph{pendent} vertex if $d(v)=1$, is \emph{isolated} if $d(v)=0$, and is \emph{dominating} if $d(v)=\nu(G)-1$. Denote by $dis(u,v)$ the distance between $u$ and $v$.\vspace*{2mm}

The \emph{complete bipartite graph} with bipartition $(X,Y)$ is denoted by $K_{|X|,|Y|}$. The complete bipartite graph $K_{1,n-1}$ is sometimes called a \emph{star } of order $n$. In addition, denote by $K_n$, $C_n$ and $P_n$ the \emph{complete} graphs, the \emph{cycle} and the \emph{path} of order $n$, respectively.  \vspace*{2mm}

The \emph{union} of  graphs $G_1=(V(G_1),E(G_1))$ and
  $G_2=(V(G_2),E(G_2))$, denoted by $G_1 \cup
G_2$, is the graph with vertex-set $V (G_1) \cup V (G_2)$ and
edge-set $E(G_1) \cup E(G_2)$.  
A graph $G_1$  is called the spanning subgraph of $G$ if $V(G_1)=V(G)$ and $E(G_1)\subseteq E(G)$. Let $V_1\subseteq V$. Denoted by $G[V_1]$ the subgraph induced by the vertex set $V_1$ and by $G\backslash V_1$ the graph obtained from $G$ by deleting $V_1$ and all edges incident to them. We sometimes write $G\backslash V_1$ as $G-v$ if $V_1$ contains exactly one element $v$.\vspace{2mm}


The following two lemmas play an important role in the following discussion.\vspace{2mm}

\begin{lemma} \label{23} {\em \cite{kr}} A graph $G$ is difference if and only if $G$ is bipartite and the neighborhoods of vertices in one of the partite sets can be linearly ordered by inclusion.
\end{lemma}

\begin{lemma} \label{24} {\em \cite[Proposition 2.5(2)]{hps}} A bipartite graph $G$ is difference if and only if $G$  contains no induced subgraphs
$P_5$.
\end{lemma}


A threshold graph can be obtained through an iterative process which
starts with an isolated vertex, and where, at each step, either a new isolated vertex is
added, or a dominating vertex is added. Then a connected threshold graph $G $ can be represented as the vector  $(0^{h_1}, 1^{h_2}, \ldots, 0^{h_{l-1}}, 1^{h_l} )$, where $l$ is even, each $h_i$ is
a positive integer number and $\sum_{i=1}^lh_i=\nu(G)$; see \cite{mp}. For convenience,   $(0^{h_1}, 1^{h_2}, \ldots, 0^{h_{l-1}}, 1^{h_l} )$ is said to be the vertex-eigenvector of the threshold graph $G$. \vspace*{2mm}

Let $G$ be a  difference graph with bipartition $(X; Y )$. Suppose that $X=\cup_{i=1}^kX_i$ and $Y=\cup_{i=1}^pY_i$ such that, for each $i$, both $X_i$ and $Y_i$ are non-empty,  all elements in $X_i$ (resp. $Y_i$) are duplicate,
  $$N(X_1)\supset N(X_2)\supset \ldots \supset N(X_k)~~and  ~~ N(Y_1)\supset N(Y_2)\supset \ldots \supset N(Y_p).$$
Then $(X_1,X_2,\ldots,X_k; Y_1,Y_2,\ldots,Y_k)$ and $(x_1,x_2,\ldots,x_k; y_1,y_2,\ldots,y_k)$ are defined as the vertex bipartition and the vertex-eigenvector of  $G$, respectively, where $|X_i|=x_i$ and $|Y_i|=y_i$ for $i=1,2,\ldots,k.$ The integer number $k$  is called  the character of  $G$; see \cite{gs}. \vspace{2mm}

For $n\ge 6$ and $n-1\le m\le 2(n-2)$, the unique $4$-Sachs minimal graph in  $\mathfrak{B}_{n,m}$ is determined in \cite{gs} as follows.
 \begin{lemma} \label{2111} \em{\cite[Theorems $23$]{gs}} Let $n\ge 6$ and $n-1<m< 2(n-2)$. Then the unique $4$-Sachs minimal graph in  $\mathfrak{B}_{n,m}$  is the  difference graph with vertex eigenvector $(1,1;m-n-2,2n-4-m)$. Moreover, the corresponding minimal $4$-Sachs number is $$\bar{\mathbf{a}}_4(\mathfrak{B}_{n,m})=(2n-4-m)(m-n+2).\eqno{(2.1)}$$
\end{lemma}\vspace{2mm}

In addition, as preliminary, we need to introduce a result on graphs having maximum number of pairs of different edges.
{\it Which graphs have
maximum number of pairs of different edges that have a common vertex among all graphs in $\mathfrak{G}_{n,m}$?} This question was first posed by Ahlswede and Katona \cite{ak} in 1978, which was completely solved in 2009 by Abrego et al. \cite{afnw}. Here we only mention a roughly characterization on those graphs  as follows:
\begin{lemma}\label{21000} {\em \cite{afnw}} If $G$ is the graph  having
maximum number of pairs of different edges that have a common vertex in $\mathfrak{G}_{n,m}$, then $G$ is a threshold graph.
\end{lemma}

\section{ A characterization on  $4$-Sachs minimal graphs}

In this section, we will give a structural characterization on the $4$-Sachs minimal graphs, showing that each $4$-Sachs minimal graph contains a difference graph as its spanning subgraph.
 Let $u$ and $v$ be two vertices of the graph $G$. Define
 $$N_G(u,v) = \{x\in V (G) \backslash \{u,v\} : xu\in E(G),xv\in E(G)\}$$ and
$$N_G(u,\bar{v}) = \{x\in V (G) \backslash \{u,v\} : xu\in E(G),xv\notin E(G)\}.$$
Let $G_{u\rightarrow v}$ be the graph formed by deleting all edges between $u$ and $N_G(u,\bar{v})$ and adding all edges from $v$ to $N_G(u,\bar{v})$. This operation is called the \emph{compression} of $G$ from $u$ to $v$; see for example \cite{kr}. It is clear that  $G_{u\rightarrow v}$ has the same number of edges as that of $G$.\vspace{2mm}

Let  $G$ be a graph with $u,v\in V(G)$.  If $dis(u,v)=2$,
the following  result tell us the compression of $G$ from $u$ to $v$ can minimize its $4$-Sachs number.\vspace{2mm}
 \begin{lemma}\em{\cite[Theorem 9]{gs}} \label{26} Let $G$ be a graph with $u,v\in V(G)$.
If $dis(u,v)=2$, then
 $$\mathbf{a}_4(G)\ge \mathbf{a}_4(G_{u\rightarrow v})$$
 inequality holds if and only if $N_G(\bar{u},v)\neq \emptyset$ and $N_G(u,\bar{v})\neq \emptyset.$
\end{lemma}\vspace{2mm}

In generally, Lemma \ref{26} does not hold for vertices $u,v$ if $dis(u,v)=1$. For  $dis(u,v)=1$, we can partly compare $\mathbf{a}_4(G)$ and $\mathbf{a}_4(G_{u\rightarrow v})$ as follows. Let $$\hat{E}_{u,v} = \{xy: xy\in E(G), x \in N_G(u,\bar{v}),y \in N_G(\bar{u},v)\} .$$

%

 \begin{theorem} \label{27} Let $G$ be a graph with $u,v\in V(G)$.
If $uv\in E(G)$ and $\hat{E}_{u,v}=\emptyset$, then
 $$\mathbf{a}_4(G)\ge \mathbf{a}_4(G_{u  \rightarrow v})$$
 inequality holds if $N_G(\bar{u},v)\neq \emptyset$ and $N_G(u,\bar{v})\neq \emptyset.$
\end{theorem}
\noindent {\bf Proof.} Let $H:= G_{u\rightarrow v}$ and $Q(G)$  denotes the set of all  quadrangles contained in  $G$. Similar to the proof of Lemma \ref{26}; see \cite[Theorem 9]{gs}, it is sufficiency to prove that  $$q(H)\ge q(G),\eqno{(3.1)}$$ where $q(\cdot)=|Q(\cdot)|$.\vspace{2mm}



To prove (3.1), we   construct an injection from $Q(G)\backslash Q(H)$ to $ Q(H)\backslash Q(G )$ in which preserves
the number of quadrangles.
The replacement function $r : E(G) \mapsto E(H)$  is defined as $$r(e)=\left \{ \begin{array}{ll}
va, & {\rm if \mbox{ }}  e=ua\mbox{ } {\rm with} \mbox{ } a\in N_G(u);\\
ub, & {\rm if \mbox{ }}  e=vb\mbox{ } {\rm with} \mbox{ } b\in N_G(u,v);\\
e, & {\rm otherwise. \mbox{ }}
\end{array}\right.$$  One can verify that $r(e)$ is an edge in $H$ for any $e\in E(G)$. 
The injection $\phi: Q(G)\backslash Q(H)\mapsto Q(H)\backslash Q(G)$   is defined by
$$\phi(C)=\{r(e):e\in C, C\in Q(G)\backslash Q(H) \}.$$
Let $C$ be a $4$-cycle of $Q(G)\backslash Q(H)$. Then $C$ must
  contain an edge $uw$ with $w\in N_G(u,\bar{v})$ and  another edge   $ux$ with $x\in N_G(u)$, that is, $C=uwyxu$ with $y\in N_G(w,x)$. By the definition of $r(e)$, $r(uw)=vw$, $r(ux)=vx$ and $r(e)=e$ if $e\notin \{uw,uy\}$. Then $\phi(C)=vwyxv$ and thus $\phi(C)\in Q(H)\backslash Q(G )$.\vspace{2mm}

Then we need to show that $\phi$ has a left inverse. Consider $r': E(H)\rightarrow E(G)$ defined by
$$r'(e)=\left \{ \begin{array}{ll}
ua, & {\rm if \mbox{ }}  e=va\mbox{ } {\rm with} \mbox{ } a\in N_G(u);\\
vb, & {\rm if \mbox{ }}  e=ub\mbox{ } {\rm with} \mbox{ } b\in N_G(u,v);\\
e, & {\rm otherwise. \mbox{ }}
\end{array}\right.$$
Define $\phi' : Q(H)\backslash Q(G)\rightarrow Q(G)\backslash Q(H)$ by $\phi'(C) = \{r'(e) : e \in C, C\in Q(H)\backslash Q(G)\}$. Then one can verify that $$\phi'(\phi(M)) = M.$$ Thus $\phi$ has a left inverse and so $\phi$ is injective. Consequently, the result follows.\hfill $\blacksquare$\vspace{2mm}

\noindent {\bf Remark 1.} The condition $\hat{E}_{u,v}=\emptyset$ in Theorem \ref{27} is necessary. Let $G$ be the graph obtained from the cycle $C_5=u_1e_1u_2e_2u_3e_3u_4e_4u_5e_5u_1$ by adding the edge  $u_1u_3$. Then $\hat{E}_{u_1,u_3}=\{e_4\}\neq \emptyset$. However, one can verify that $a_4(G)=4$ and $a_4(G_{u_1\rightarrow u_3})=5.$ \vspace{2mm}


Combining with Lemmas  \ref{24}, \ref{26} and Theorem \ref{27}, we  give a structural characterization on  $4$-Sachs minimal graphs in $\mathfrak{G}_{n,m}$ as follows.\vspace{2mm}

 \begin{theorem} \label{12} Let $G$ be a $4$-Sachs minimal graph in $\mathfrak{G}_{n,m}$. Then $G$ contains a difference graph as its spanning subgraph.
\end{theorem}
\noindent {\bf Proof.} From Lemma \ref{24}, $G$ contains no induced subgraph $P_5$, then by Lemma \ref{26} $G$ itself is difference if $G$ is bipartite.
On the other hand, the result follows  if $G$ is complete as each complete graph contains the star $K_{1,n-1}$, a difference graph, as its spanning subgraph. \vspace{2mm}

Suppose now that $G$ is  non-bipartite and non-complete.
 Let $u_1$ be an arbitrary vertex with $d_G(u_1)<n-1$ and let $S=\{u_i|i=1,2,\ldots,s\}$ be the maximal independent set
 such that for each $i(i\neq 1)$ the distance between $u_1$ and $u_i$ is even. (Since $G$ is non-bipartite and non-complete, such a vertex $u_1$ must exist and $S$ contains at least two elements, $u_1$ itself and  another vertex.)\vspace{2mm}

 We  claim that $dis(u_i,u_j)=2(i\neq j)$.
 By Lemma \ref{24}, $G$ contains no induced subgraphs $P_5$, then $dis(u_1,u_i)=2$ for $i=2,3,\ldots,s$ and, for any pair vertices $u_i$ and $u_j$, $dis(u_i,u_j)\in\{2,3\}$. Assume  that there exists vertices $u_i$ and $u_j$ such that $dis(u_i,u_j)=3$. Then, applying Lemma \ref{24} again, there exist vertices $v$ and $w$ such that $v\in N_G(u_1,u_i)$, $w\in N_G(u_1,u_j)$  and $u_ivwu_j$ forms a distance path between $u_i$ and $u_j$.
 Moreover, we have  $\hat{E}_{v,w}=\emptyset$.
 (Otherwise, assume to the contrary that $x\in N_G(v,\bar{w})$ and $y\in N_G(w,\bar{v})$ such that $xy\in E(G)$, then $dis(x,u_1)=2$ with $N_G(x)\nsubseteq N_G(u_1)$ and $N_G(u_1)\nsubseteq N_G(x)$, a contradiction.) Thus by Theorem \ref{27} $a_4(G_{v\rightarrow w})<a_4(G)$, a contradiction is yielded.   Consequently, the claim is true. \vspace{2mm}

  Applying Lemma \ref{26},  the neighborhoods of vertices in $S$ can be linearly ordered by inclusion.  Without loss of generality, suppose that
$$N_G(u_1)\subseteq N_G(u_2) \subseteq \ldots  \subseteq N_G(u_{s}).\eqno{(3.2)}$$

%
%

If $G[V\backslash S]$ is  complete, then $G$ is threshold and thus the result follows from the fact that an arbitrary threshold graph can be obtained from a difference graph
by adding all possible edges in one of the partite sets. If $G[V\backslash S]$ is  non-complete. By the method similar to  above, we get another maximal independent set $T=\{v_i|i=1,2,\ldots,t\}$  of $ V\backslash S$ and  the distance between arbitrary two distinct vertices   of $T$ is $2$.
  Applying Lemma \ref{26} again,  the neighborhoods of vertices in $T$ can be linearly ordered by inclusion. Suppose that
$$N_G(v_1)\subseteq N_G(v_2) \subseteq \ldots  \subseteq N_G(v_{t}).\eqno{(3.3)}$$


Then the graph $G[S\cup T]$ is a bipartite graph with bipartition $(S; T)$.  Otherwise $G[S\cup T]$ contains triangles, which  contradicts to that both $S$ and $T$ are independent sets. Consequently, combining with (3.2) and (3.3) $G[S\cup T]$ is a difference graph.\vspace{2mm}

Recall that $G$ is non-bipartite, then $V\backslash (S\cup T)\neq \emptyset$, say $x\in V\backslash (S\cup T)$. Note that $x\notin S$ and $x\notin T$, then $x$ is adjacent to some vertex of $S$ and some vertex of $T$. Thus $x\in N_G(u_s)$ by Eq. (3.2) and $x\in N_G(v_t)$ by Eq. (3.3). Assume that $N_G(x)\cap S\neq \emptyset$, then
$u_1\notin N_G(x)$. Applying Theorem \ref{27}  $N_G(u_1)\subseteq N_G(x)$ as $u_s\in N_G(x)$ and thus we can delete some adjacent edges of $x$ such that $S\cup T \cup \{x\}$ forms a difference graph. By a similar method, the result follows if $V\backslash (S\cup T\cup \{x\})\neq \emptyset$. Consequently, the proof is complete.\hfill $\blacksquare$\vspace{2mm}

\noindent{\bf Remark 2.} For a given $4$-Sachs minimal  $(n,m)$-graph, its spanning difference subgraph may not unique. Let $G$ be the threshold graph with vertex-eigenvector $(0^2,1^2)$,  $K_4-e$. Then $G$ is the graph having minimal  $4$-Sachs number among all graphs in  $\mathfrak{G}_{4,5}$. One can verify that each of the difference graphs $K_{2,2}$ and $K_{1,3}$ is a spanning difference subgraph of $G$. \vspace{2mm}

\section{  $4$-Sachs minimal  graphs for $n\ge 6$ and $n-1\le m\le 2n-4$.}
For $n\ge 6$ and $n-1\le m\le 2n-4$, we in this section  determine all $4$-Sachs minimal  $(n,m)$-graphs  together with the corresponding minimal $4$-Sachs number. \vspace{2mm}

Let $G$ be a $4$-Sachs minimal  $(n,m)$-graph. Hereafter, we always use the notation  $\mathbf{B}_G$ to denote the given spanning difference subgraph of $G$ and use $\overline{\mathbf{B}_G}$ to denote the  graph obtained from $G$ by deleting all edges of $\mathbf{B}_G$ together with all isolated vertices of the resulting graph. Then $G= \mathbf{B}_G\cup \overline{\mathbf{B}_G}$ and $E(\mathbf{B}_G)\cap E(\overline{\mathbf{B}_G})=\emptyset$. We begin our discussion with a lower bound on $\epsilon(\overline{\mathbf{B}_G})$, with respect to the given spanning difference subgraph $\mathbf{B}_G$.\vspace{2mm}

 \begin{lemma} \label{281}  Let $G$ be a $4$-Sachs minimal  $(n,m)$-graph and $\mathbf{B}_G$ be a spanning difference subgraph of $G$ with  vertex bipartition  $(X_1,X_2,\ldots,X_k;Y_1,Y_2,\ldots,Y_k) (k\ge 2)$. Let also the graph $\overline{\mathbf{B}_G}$ is defined as above. If there exists an integer $p(2\le p\le k)$ such that $t_p:=|V(\overline{\mathbf{B}_G})\cap X_p|>0$, then
$$\epsilon(\overline{\mathbf{B}_G})\ge \sum_{j=1}^{p-1}|X_j|+t_p-1.$$
  \end{lemma}
\noindent {\bf Proof.} Let $x_ix_p\in \overline{\mathbf{B}_G}$ with $x_p\in X_p$ and $x_i\in X_i(i\le p)$. We claim that each vertex of $\cup_{j=1}^{p-1}X_j$ is either adjacent to $x_i$ or to $x_p$. Assume to the contrary that there exists a vertex $x\in X_t$ with $t<p$ such that $xx_i\notin  \overline{\mathbf{B}_G}$ and $ xx_p\notin  \overline{\mathbf{B}_G}$. Then neither $N_G(x)\subseteq N_G(x_p)$ nor $N_G(x_p)\subseteq N_G(x)$ and thus $\mathbf{a}_4(G_{x_p\mapsto x})<\mathbf{a}_4(G)$ by Lemma \ref{26}, a contradiction. Consequently, $E(\overline{\mathbf{B}_G}[\cup_{j=1}^{p-1}X_j\cup \{x_p\}])$ contains at least $\sum_{j=1}^{p-1}|X_j|$ edges. In addition, there has at least $t_p-1$ additional edges  incident to the remaining $t_p-1$ vertices, other than $x_p$, of $X_p$. Thus, the result follows. \hfill $\blacksquare$\vspace{2mm}

Let $G$ be a $4$-Sachs minimal  $(n,m)$-graph
with the spanning difference subgraph $\mathbf{B}_G$.  Suppose further that the bipartition of $\mathbf{B}_G$ is $(U,W)$. If $\overline{\mathbf{B}_G}$ is connected, then either $V(\overline{\mathbf{B}_G})\subseteq U$ or $V(\overline{\mathbf{B}_G})\subseteq W$ as $\overline{\mathbf{B}_G}$  contains no edges joining $U$ and $W$. Without loss of generality, suppose that $V(\overline{\mathbf{B}_G})\subseteq U$. Let now $e_i$ and $e_j$ be two adjacent edges of $\overline{\mathbf{B}_G}$, say $e_i=uu_i$ and $e_j=uu_j$. We define $$\mathbf{N}_W(e_i,e_j)=\{w|w\in W, wu_i\in G, wu_j\in G\},$$ that is, $\mathbf{N}_W(e_i,e_j)=\mathbf{N}_W(u_i,u_j)$,  the common neighbors of vertices $u_i$ and $u_j$ in the set $W$. Let $n_W(e_i,e_j)=|N_W(e_i,e_j)|.$
 Then we can give  another formula on $4$-Sachs number of a $4$-Sachs minimal $G$ as follows.

 \begin{theorem} \label{29}  Let $G$ be a $4$-Sachs minimal  $(n,m)$-graph   and $\mathbf{B}_G$ be a spanning difference subgraph of $G$ with bipartition  $(U,W)$.
 If  $V(\overline{\mathbf{B}_G})\subseteq U$, then
 $$\mathbf{a}_4(G)=\mathbf{a}_4(\mathbf{B}_G)+\mathbf{a}_4(\overline{\mathbf{B}_G})+\sum_{e_i\in \overline{\mathbf{B}_G}}\epsilon(\mathbf{B}_G\backslash V(e_i))-2\sum_{e_p,e_q\in \overline{\mathbf{B}_G}}n_W(e_p,e_q),$$ where the first summation is over all edges $e_i$ of $\overline{\mathbf{B}_G}$ and the second summation is over all adjacent edges $e_p$ and $e_q$ of $\overline{\mathbf{B}_G}$.
\end{theorem}
\noindent {\bf Proof.} We divide  all $4$-Sachs subgraphs of $G$  into four types: those that contained in $\mathbf{B}_G$; those that contained in $\overline{\mathbf{B}_G}$;  those that each of them is a $2$-matching, which contains exactly one edge of $\mathbf{B}_G$ and exactly one edge of $\overline{\mathbf{B}_G}$; and those that each of them is a quadrangle,  which contains at least one edges of $\mathbf{B}_G$ and at least one edge of $\overline{\mathbf{B}_G}$. Obviously, the $4$-Sachs number of the first type  is $\mathbf{a}_4(\mathbf{B}_G)$, the $4$-Sachs number of the second type  is $\mathbf{a}_4(\overline{\mathbf{B}_G})$, the $4$-Sachs number of the third type  is $\sum_{e_i\in \overline{\mathbf{B}_G}}\epsilon(\mathbf{B}_G\backslash V(e_i))$.
 Recall that $V(\overline{\mathbf{B}_G})\subseteq U$ by hypothesis, then each $4$-Sachs subgraph of the fourth type contains exactly two adjacent edges of $\overline{\mathbf{B}_G}$ and exactly two adjacent edges of $\mathbf{B}_G$, then  the $4$-Sachs number of those subgraphs is $-2\sum_{e_p,e_q}n_W(e_p,e_q),$ where the the summation is over all adjacent edges $e_p$ and $e_q$ of $\overline{\mathbf{B}_G}$.
  Consequently, the result follows.\hfill $\blacksquare$\vspace{2mm}

\begin{theorem} \label{29011}  Let $G$ be a graph  with $v\in V(G)$. If $d_G(v)\le 2$, then  $$\mathbf{a}_4(G)\ge \mathbf{a}_4(G-v).$$
\end{theorem}
\noindent {\bf Proof.} For $d_G(v)=1$, let $N_G(v)=\{u\}$.  From Eq.(1.1), $\mathbf{a}_4(G)=\mathbf{m}_2(G)-2\mathbf{q}(G)$, then $\mathbf{a}_4(G)= \epsilon(G- u)+ \mathbf{a}_4(G-v)\ge \mathbf{a}_4(G-v)$ as $v$ does not contained in any quadrangle of $G$. \vspace{2mm}

For $d_G(v)=2$, say $N_G(v)=\{x,y\}$, then all $2$-matchings of $G$ can be divided into the following three types: those that each of them contains the edge $vx$; those that each of them contains the edge $vy$ and those otherwise. Note that the cardinality of the former is $\epsilon(G-v)-d_{G-v}(x)$,  of the second type is $\epsilon(G-v)-d_{G-v}(y)$, and   of the third type is $\mathbf{m}_2(G-v)$, then $$\mathbf{m}_2(G)=2\epsilon(G-v)-d_{G-v}(x)-d_{G-v}(y)+\mathbf{m}_2(G-v).$$ On the other hand, we have
 \begin{equation*}
\begin{array}{lll} \mathbf{q}(G)&=&\mathbf{q}(G-v)+\mathbf{q}(G,v)\\&=&\mathbf{q}(G-v)+|N_{G-v}(x,y)|,
\end{array}
\end{equation*} where $\mathbf{q}(G,v)$ denotes the cardinality of all quadrangles, of $G$, containing the vertex $v$.
Consequently,
 \begin{equation*}
\begin{array}{lll} \mathbf{a}_4(G)&=&\mathbf{m_2}(G)-2\mathbf{q}(G)\\&= &\mathbf{m_2}(G-v)-2\mathbf{q}(G-v)+2\epsilon(G-v)-d_{G-v}(x)-d_{G-v}(y)-2|N_{G-v}(x,y)|\\&\ge &\mathbf{m_2}(G-v)-2\mathbf{q}(G-v)\\ &=& \mathbf{a}_4(G-v),
\end{array}
\end{equation*} as $\epsilon(G-v)-d_{G-v}(x)\ge |N_{G-v}(x,y)|$ and $\epsilon(G-v)-d_{G-v}(y)\ge |N_{G-v}(x,y)|$.
 Thus the result follows.\hfill $\blacksquare$\vspace{2mm}

\noindent {\bf Remark 3.} In Theorem \ref{29011}, the restriction that $d_G(v)\le 2$ is necessary. Let $v$ be an arbitrary vertex of $K_4$. Then $d_{K_4}(v)= 3$. However, we find that $$\mathbf{a}_4(K_4)=-3\le \mathbf{a}_4(K_4-v)=\mathbf{a}_4(K_3)=0.$$

 Below we focus on determining all $4$-Sachs minimal $(n,m)$-graphs with $n\ge 4$ and  $n-1\le m\le 2n-4$. Firstly, we have
\begin{theorem} \label{29012}  Let $n\ge 6$, $n-1\le m \le 2n-4$ and $G$ be an arbitrary connected $(n,m)$-graph. Then  $$\mathbf{a}_4(G)\ge 0.$$
\end{theorem}
\noindent {\bf Proof.}  Obviously, it is sufficiency to show that the inequality  holds for each  $4$-Sachs minimal graph. Thus we suppose that $G$ is a $4$-Sachs minimal $(n,m)$-graph. \vspace{2mm}

We first claim that
$\delta(G)\le 2$. By Theorem \ref{12}, $G$ contains a spanning difference subgraph. Then there exists a pair vertices, say $x$ and $y$, such that $xy\in E(G)$ and $d(x)+d(y)\ge n$. Assume that $\delta(G)\ge 3$. Thus
$$\sum_{v\in V(G)} d(v)\ge n+3(n-2)=4n-6,$$ which implies that $m\ge 2n-3$, a contradiction to the hypothesis.\vspace{2mm}

Then we prove $\mathbf{a}_4(G)\ge 0$ by induction on $n$. Obviously, the result follows if $\nu(G)=4$. Suppose that the result follows for  $\nu(G)<n(\ge 4)$. For $\nu(G)=n$, let $v\in V(G)$ such that $d(v)=\delta(G)$. Then $d(v)\le 2$ by the discussion above and thus by Theorem \ref{29011} and the inductive  hypothesis $$\mathbf{a}_4(G)\ge \mathbf{a}_4(G-v)\ge 0.$$  Consequently, the result follows. \hfill $\blacksquare$\vspace{2mm}

 \begin{theorem} \label{210}   Let $n\ge 6$,  $n-1\le m\le 2n-4$ and $G$  be a $4$-Sachs minimal graph in $\mathfrak{G}_{n,m}$.
If $\triangle(G)=n-1$, then $G$ is a threshold graph. Moreover, the vertex-eigenvector of $G$ is  $(0^{1},1^2,0^{n-4},1^1)$ or $(0^{1},1^{1},0^{n-1},1^1)$ if $m=n+2$, and  $(0^{m-n-1},1^1,0^{2n-m-3 },1^1)$ otherwise.
\end{theorem}
\noindent {\bf Proof.} Let  $d_G(v)=\triangle(G)=n-1$ and $\mathbf{B}_G:=K_{1,n-1}$, whose  central vertex be $v$.
Then $\mathbf{B}_G$ is  a spanning difference subgraph of $G$.
Let now $W=V(G)\backslash \{v\}$ and $U=\{v\}$. Then $V(\overline{\mathbf{B}_G})\subseteq W$ and $\epsilon(\overline{\mathbf{B}_G})=m-n+1\le n-3$.
  Moreover, we find that $n_U(e_p,e_q)=1$ for each pair of adjacent edges $e_p$ and $e_q$ of $\overline{\mathbf{B}_G}$, and $\epsilon(\mathbf{B}_G\backslash V(e_i))=n-3$ for each  $e_i$ of $\overline{\mathbf{B}_G}$ (if there exists),
   then to minimize $\mathbf{a_4}(G)$, we need to minimize $\mathbf{a_4}(\overline{\mathbf{B}_G})$ and to maximize the number of  pairs of adjacent edges contained in $\overline{\mathbf{B}_G}$. Recall that $\epsilon(\overline{\mathbf{B}_G})\le n-3$, then from Lemma \ref{21000} the graph $\overline{\mathbf{B}_G}$ can be chosen as a threshold graph such that whose $4$-Sachs number is zero as $|W|=n-1$. Thus $\overline{\mathbf{B}_G}$ is $C_3$ if $\overline{\mathbf{B}_G}\in \mathfrak{G}_{3,3}$ and $\overline{\mathbf{B}_G}=K_{1,m-n+1}$ otherwise. Consequently, the result follows.\hfill $\blacksquare$\vspace{2mm}

 Let $G$  be a $4$-Sachs minimal graph in $\mathfrak{G}_{n,m}$ with $n\ge 6$,  $n-1\le m\le 2n-4$ and $\triangle(G)=n-1$. Then from Theorem \ref{210} and by a directly calculation, we have
  $$\mathbf{a}_4(G)=(m-n+1)(2n-m-3).\eqno{(4.3)}$$

\begin{lemma}\label{21001}  Let $n\ge 6$ and $G$  be a $4$-Sachs minimal  graph in $\mathfrak{G}(n,2n-4)$. If $\Delta(G)= n-2$, then $G$ is the complete bipartite graph   $K_{2,n-2}$.
\end{lemma}
\noindent {\bf Proof.} By Lemma \ref{2111} the result follows if $G$ is bipartite. Then we need only to show that $G$ is indeed bipartite. Assume to the contrary that $G$ is non-bipartite. Combining with  Theorem \ref{12} and the fact that $\Delta(G)= n-2$, then $G$ contains the difference graph with vertex-eigenvector $(1,1;p,n-2-p)$, denoted by $\mathbf{B}_G$, as its spanning subgraph with $1\le p\le n-3$. Further, we assume that the vertex bipartition of $\mathbf{B}_G$ is $(U_1,U_2; W_1,W_2)$ with $U_1=\{u_1\}$, $U_2=\{u_2\}$, $W_1=\{w_i|i=1,2,\ldots,p\}$ and $W_1=\{w_i|i=p+1,\ldots,n-2\}$.   Then one of the following two cases must occur:\\
 {\bf Case 1.} $G$ contains no vertices with degree two. \vspace{2mm}

\noindent Then the degree of each vertex of $G$ is either $1$ or at least $3$. Thus $W_1\subseteq V(\overline{\mathbf{B}_G})$, where $\overline{\mathbf{B}_G}$ is defined as above. By Lemma \ref{21000} $\overline{\mathbf{B}_G}$ is a threshold graph. Thus there exists a vertex, say $w_1$, such that $d_{\overline{\mathbf{B}_G}}(w_1)=\nu(\overline{\mathbf{B}_G})-1$ with $\nu(\overline{\mathbf{B}_G})\ge p$. Denote by $G'$ the graph obtained from $G$ by deleting all edges between $u_1$ and $W_1 \backslash V(\overline{\mathbf{B}_G})$ and adding all edges between $w_1$ and $W_1 \backslash V(\overline{\mathbf{B}_G})$. Then $\mathbf{a}_4(G')<\mathbf{a}_4(G)$ by Lemma \ref{27} and $\Delta(G')=n-1$. From Eq. (4.3) $\mathbf{a}_4(G')\ge n-3>0$, which is contradiction to the assumption. Thus such a case can not be occur.\\
 {\bf Case 2.} $G$ contains vertices with degree two. \vspace{2mm}

 \noindent Let $v\in V(G)$ be the vertex with degree $2$, say $N_G(v)=\{u_1,u_2\}.$ Denote by $q(G,v)$ the number of quadrangles containing the vertex $v$. Then $q(G,v)=|N_{G-v}(u_1,u_2)|$ and the number of $2$-matchings, of $G$, containing the edge $vu_1$ is $m-d_{G-v}(u_1)-2$ and the number of $2$-matchings, of $G$, containing the edge $vu_2$ is $m-d_{G-v}(u_2)-2$.
Consequently, $$\mathbf{a}_4(G)=2m-d_{G-v}(u_1)-d_{G-v}(u_2)-4-2|N_{G-v}(u_1,u_2)|+\mathbf{a}_4(G-v).$$
Note that $\epsilon(G-v)=2n-6=2(\nu(G-v)-2)$ and $\Delta(G)= n-2$, then by Theorem \ref{29012} $\mathbf{a}_4(G-v)\ge 0$ and by the minimality of $\mathbf{a}_4(G)$, $$d_{G-v}(u_1)=d_{G-v}(u_2)=|N_{G-v}(u_1,u_2)|=n-3.$$Thus, $G=K_{2,n-2}$. Consequently, the result follows.\hfill $\blacksquare$\vspace{2mm}

  \begin{theorem} \label{211}  Let $n\ge 6$,  $n-1\le m\le 2n-4$ and $G$  be a $4$-Sachs minimal  graph in $\mathfrak{G}_{n,m}$. If $\Delta(G)< n-1$, then $G$ is the difference graph with vertex eigenvector $(1,1;m-n-2,2n-4-m)$.
\end{theorem}
\noindent {\bf Proof.} Let  $d(u_1)=\Delta(G):=p.$ 
  Suppose that $V=U\cup W$ with $W=N(u_1)=\{w_i|i=1,2,\ldots,p\}$ and $U=V\backslash W=\{u_i|i=1,2,\ldots,n-p\}.$\vspace{2mm}

   We first show that the graph, denoted by $\mathbf{B}$, obtained from $G$ by deleting all edges of $G[W]$ (if there exists) is a spanning difference subgraph of $G$.
By the definition of the set $U$, $dis(u_1,u_i)\ge 2$ for $i=2,\ldots,n-p$. If there exists a vertex $u_j$ with $dis(u_1,u_j)> 2$, then there exists a vertex $u_i$ such that $dis(u_1,u_i)= 2$ and $N(u_i)\backslash N(u_1)\neq \emptyset$. Then $N(u_1)\backslash N(u_i)\neq \emptyset$ as $d(u_i)\le \Delta(G)=d(u_1).$
Thus applying Lemma \ref{26} the graph $G_{u_1\rightarrow u_i}$ has less $4$-Sachs number than that of $G$ and $G_{u_1\rightarrow u_i}\in \mathfrak{G}_{n,m}$, which is a contradiction. Consequently,  $dis(u_1,u_j)= 2$ for each $j(2\le j\le n-p)$. Furthermore,   the neighborhoods of any vertex in $U$ are linearly ordered by inclusion. On the contrary, assume that there exist vertices $u_i$ and $u_j$ such that
$w_i\in N(u_i), w_j\in N(u_j), w_i\notin N(u_j)$ and $w_j\notin N(u_i)$. Then, no matter $w_iw_j\in E(G)$ or not, a graph having less $4$-Sachs number can be obtained applying Lemma \ref{26} or Theorem \ref{27}, which is also a contradiction.
Without loss of generality, suppose that $$W=N_G(u_1)\supseteq N_G(u_2) \supseteq \ldots  \supseteq N_G(u_{n-p}).\eqno{(4.4)}$$
 (4.4) compels the neighborhoods of vertices of $W$ in $\mathbf{B}$ are linearly ordered by inclusion.
Thus  $\mathbf{B}$ is a difference graph.  Without loss of generality, suppose that $$U=N_B(w_1)\supseteq N_B(w_2) \supseteq \ldots  \supseteq N_B(w_{p}).\eqno{(4.5)}$$

Recall that $n-1\le m \le 2n-4$, then $G$  is the difference graph with vertex eigenvector $(1,1;m-n-2,2n-4-m)$ if $E(\overline{\mathbf{B}})=\emptyset$ by Lemma \ref{2111}. Thus it remain to show that $E(\overline{\mathbf{B}})$ is indeed an empty set. \vspace{2mm}

Assume that the vertex bipartition of $\mathbf{B}$ is $(U_1,\ldots,U_k;W_1,\ldots,W_k),$  where $\cup_{i=1}^kU_i =U$ and $\cup_{i=1}^kW_i =W$.
  For each $i(i=1,2,\ldots,k)$, suppose that $U_i=\{u_i^{s_i}|s_i=1,\ldots,|U_i|\}$ and $W_i=\{w_i^{t_i}|t_i=1,\ldots,|W_i|\}$. Since $U$ is independent, $N(u_1^{1})= W$ and $U\subseteq N(w_1^{1})$. Then $d_G(u_1^{1})=\ldots=d_G(u_1^{|U_1|})=\Delta(G)=p$. Applying Theorem \ref{27}, $d_G(u_k^1)=\cdots=d_G(u_k^{|U_k|})\ge 2$ and thus $|W_1|\ge 2$. Otherwise, assume that $d_G(u_k^1)=\cdots=d_G(u_k^{|U_k|})=1$, that is, each vertex of $U_k$ is a pendent vertex of $G$, then the graph obtained from $G$ by deleting all edges between $w_1^1$ and $U_k$ and adding all edges between $u_1^1$ and $U_k$ has less $4$-Sachs number.
  Thus $\mathbf{B}$ contains a difference graph $\mathbf{B}'$ with vertex eigenvector $(1,1,|U|-2;2,\sum_{i=1}^{k-1}|W_i|-2,|W_k|)$ as its spanning subgraph. Let $e(\overline{\mathbf{B}})=t$. Then $\mathbf{B}'\in \mathfrak{B}_{n,m'}$ with $m'\le m-t$, and thus $$m'=2\sum_{i=1}^k|W_i|-|W_k|+2|U|-4\le m-t.$$ Recall that $\sum_{i=1}^k|W_i|+|U|=|W|+|U|=n$, then
  $$|W_k|\ge 2n-4-m+t,$$ which implies that $\mathbf{B'}$, as well as $\mathbf{B}$, contains at least $2n-4-m+t$ pendent vertices.

  We further claim that $G$  contains at least $2n-4-m+t$ pendent vertices. The result follows if each pendent vertex of $\mathbf{B}$ is a pendent vertex of $G$. Assume that $l$ pendent vertices of $\mathbf{B}$ are not still pendent vertices of $G$. By Lemma \ref{281} $t\ge \sum_{i=1}^{k-1}|W_i|+l-1$, then $t$ increase with $l$ and thus this operation will yield more other  pendent vertices. Consequently, the claim is true.\vspace{2mm}

By the assumption $t\ge 0$, then  $G$  contains at least $2n-4-m$ pendent vertices. Moreover, from Lemma \ref{26}, all those pendent vertices have the same neighbor $u_1^1$.
 Denote by $G'$ the graph obtained from $G$ by deleting $2n-4-m$ pendent vertices together all edges incident to them. Then $\nu(G')=m-n+4$, $\epsilon(G')=2m-2n+4$ and thus $\epsilon(G')= 2\nu(G')-4$. Consequently, $\mathbf{a_4}(G')\ge 0$ by Theorem \ref{29012}.
Furthermore,
$$\mathbf{a}_4(G)= (2n-4-m)(m-d_G(u_1^1))+\mathbf{a}_4(G').$$ Recall that $d_G(u_1^1)\le \triangle(G)<n-1$, then, to minimizes $\mathbf{a_4}(G)$, $d_G(u_1^1)=n-2$ and $\mathbf{a_4}(G')=0$ by Theorem \ref{29012}.  Applying Lemma \ref{21001}, $G'$ is the complete bipartite graph $K_{2,m+4-n}$, which imples that $t=0$. Consequently, the proof is completed. \hfill $\blacksquare$\vspace{2mm}

Combining with Theorems \ref{210} and \ref{211}, we determine all connected $4$-Sachs minimal $(n,m)$-graph as follows.
\begin{theorem} \label{2100}  Let   $n\ge 6$ and  $n-1\le m\le 2n-4$. Let  $G$ be the $4$-Sachs minimal $(n,m)$-graph. Then   $$\mathbf{\bar{a}_4}(\mathfrak{G}_{n,m})=\left \{ \begin{array}{ll}
(2n-3-m)(m-n+1), & {\rm if  \mbox{ } n\le m\le \frac{3n-5}{2}};\\
(2n-4-m)(m-n+2), & {\rm if  \mbox{ } \frac{3n-5}{2}\le m\le 2n-4.}
 \end{array}\right.$$
Moreover, $$G= \left \{ \begin{array}{ll}
G_1~or~G_2, & {\rm if  \mbox{ }  m=n+2};\\
G_2, & {\rm if  \mbox{ } n\le m< \frac{3n-5}{2}\mbox{ } and \mbox{ }m\neq n+2};\\
G_2~or~G_3, & {\rm if  \mbox{ }  m=\frac{3n-5}{2}};\\
G_3, & {\rm if  \mbox{ } \frac{3n-5}{2}< m\le 2n-4,}
 \end{array}\right.$$
where $G_1$ and $G_2$ are threshold graphs with vertex eigenvector $(0^1,1^2,0^{n-4},1^1)$ and $(0^{m-n-1},1^1,0^{2n-m-3},1^1)$, respectively, and $G_3$ is the difference graph with vertex eigenvector $(1,1;m-n-2,2n-4-m)$.
\end{theorem}
\noindent {\bf Proof.} If $\Delta(G)=n-1$, then by Theorem \ref{210} $G$  is the threshold graph with vertex-eigenvector   $(0^{1},1^2,0^{n-4},1^1)$ or $(0^{1},1^{1},0^{n-1},1^1)$ if $m=n+2$, and  $(0^{m-n-1},1^1,0^{2n-m-3 },1^1)$ otherwise. If $\Delta(G)<n-1$, then by Theorem \ref{211} $G$  is the difference graph with vertex eigenvector $(1,1;m-n-2,2n-4-m)$. Thus the result follows by comparing Eq.s(2.1) and (4.3). \hfill $\blacksquare$\vspace{2mm}

\end{document}